\documentclass[11pt,a4paper,mathscr,twoside]{amsart}

\usepackage{amsthm}
\usepackage{amsfonts, amssymb, amsmath, amscd, latexsym}
\usepackage{psfrag, eucal, enumerate, latexsym}
\usepackage[all]{xy}

\newtheorem{thm}{Theorem}[section]

\newtheorem{lem}[thm]{Lemma}
\newtheorem{corol}[thm]{Corollary}
\newtheorem{prop}[thm]{Proposition}

\theoremstyle{definition}
\newtheorem{rmk}[thm]{Remark}

\newtheorem{dfn}[thm]{Definition}

\newcommand{\Z}{\mathbb Z} \newcommand{\C}{\mathbb C} \newcommand{\Q}{\mathbb Q}
 
\newcommand{\p}{\mathbb P} \newcommand{\pd}{\check{\mathbb P}}
\newcommand{\G}{\mathbb G}
 
  \newcommand{\kk}{k} 

\newcommand{\OO}{\mathcal{O}}  
  \newcommand{\KK}{\mathcal{K}}
\newcommand{\CC}{\mathcal{C}}  \newcommand{\DD}{\mathcal{D}}

\DeclareMathOperator{\rk}{rk}

\DeclareMathOperator{\Hilb}{Hilb} 

\DeclareMathOperator{\ts}{\otimes} \DeclareMathOperator{\ev}{ev}
 
 \DeclareMathOperator{\Lm}{L} \DeclareMathOperator{\Rm}{R} 
\DeclareMathOperator{\D}{D} 
 
\DeclareMathOperator{\VPS}{VPS}
\DeclareMathOperator{\supp}{supp}

\DeclareMathOperator{\SL}{{\sf SL}}

\DeclareMathOperator{\SO}{{\sf SO}}

\DeclareMathOperator{\PGL}{{\sf PGL}}

\DeclareMathOperator{\Ll}{{\mathbf L}}
\DeclareMathOperator{\Rr}{{\mathbf R}}

\DeclareMathOperator{\Ext}{Ext} \DeclareMathOperator{\Tor}{Tor} \DeclareMathOperator{\Sym}{Sym} 
\DeclareMathOperator{\Hom}{Hom} \DeclareMathOperator{\im}{Im} \DeclareMathOperator{\cok}{coker} \DeclareMathOperator{\End}{End}
 
\DeclareMathOperator{\HH}{H}
\DeclareMathOperator{\hh}{h}  
 \DeclareMathOperator{\M}{M}

\DeclareMathOperator{\V}{V}  

\newcommand{\rt}{\longrightarrow}       
       \newcommand{\sr}{\rightarrow}
     
\newcommand{\wq}{\wedge^2 Q^*}
                      
\newcommand{\bx}{\boxtimes}

\SelectTips{cm}{}
\newdir{ >}{{}*!/-5pt/@{>}}

\begin{document}


\title{Bundles over the Fano threefold $V_5$}

\author{Daniele Faenzi}
\address{Dipartimento di Matematica ``U.~Dini'', Universit\`a di
 Firenze, Viale Morgagni 67/a, I--50134, Florence, Italy}
\email{faenzi@math.unifi.it}
\urladdr{http://www.math.unifi.it/\~{}faenzi/}



\begin{abstract}
Using an explicit resolution of the diagonal for the variety $V_5$,
we provide cohomological characterizations of the universal and quotient bundles.
A splitting criterion for bundles over $V_5$ is also proved.

The presentation of semistable aCM bundles is shown, 
together with a resolution--theoretic classification of low rank aCM bundles.
\end{abstract}



\sloppy

\maketitle

\section{Introduction}

According to the classification of Fano threefolds with second Betti
number equal to $1$ (see \cite{iskovskih:II} and \cite{iskovskih:I}), 
there are $4$ Fano threefolds with $b_2 =1$ and $b_3=0$. Namely the
projective space $\p^3$, the quadric $Q_3$, the degree-$5$ {\em del
Pezzo} threefold $V_5$ and the genus-$12$ variety $V_{22}$. The
cohomology ring over $\C$ of any of these is isomorphic to the
cohomology ring of $\p^3$ i.e. $\C[h] / h ^4 = \oplus
\HH^{p,p}(X)$, where $h$ is the hyperplane divisor, so  we say
that these are the cases with {\em trivial} Hodge structure.

The cohomology rings are nonisomorphic over $\Z$ by the degree of 
$h^2$ with respect to the generator of $\HH^4(X)$, the Poincar\'e
dual of a straight line. This degree is the degree of the threefold:
respectively $1$, $2$, $5$, $22$. 
The intermediate Jacobian is trivial and the Chow ring is isomorphic
to the cohomology ring.

All of these varieties are rational. The variety $V_{22}$ is the only one with nontrivial
moduli space and this moduli space admits a {\em special point} corresponding
to a distinguished threefold $U_{22}$. The variety $U_{22}$ is called the
Mukai--Umemura threefold.

Any of the varieties $\p^3$, $Q_5$, $V_5$ and $U_{22}$ admits a quasi--homogeneous
structure for the group $\SL(2)$ (equivalently for $\PGL(2)$ or
$\SO(3)$), i.e. it is the (smooth) closure of an orbit by the action of
this group. We will briefly resketch this in section (\ref{quasi-homogeneous}).

However here we will concentrate on $V_5$ and on vector
bundles over it. In section (\ref{bundles}) we review some facts about
bundles over $V_5$, while in section (\ref{resolution}) we will give
an explicit resolution of the diagonal, see Theorem \ref{V5-main}.
This is constructed rather directly, but, unlike the case of projective spaces and
Grassmannians, cfr. \cite{beilinson:derived-and-linear} and 
\cite{kapranov:derived-homogeneous}, it is not given by a Koszul complex. Results 
for $V_5$ are obtained also by Orlov in \cite{orlov:V5}, while Canonaco in
\cite{canonaco:beilinson-weighted} provides similar results for weighted projective spaces 
and Kawamata in \cite{kawamata:derived-stacks} exhibits
an infinite resolution for general projective varieties.

In section (\ref{splitting-section}) we prove a splitting criterion making use of
the canonical resolution, and cohomological characterizations of
universal and quotient bundles. Then we concentrate on bundles with the
property that $\HH^p(V_5,E \ts \OO(t))=0$ for any $t \in \Z$ and $0<p<3 = \dim(V_5)$. We
call a vector bundle $E$ with this property {\em aCM} for arithmetically Cohen--Macaulay,
following standard terminology.

We describe moduli spaces of aCM bundles over $V_5$ using resolutions and prove that 
semistable aCM bundles admit a peculiar presentation that may lead to classification 
at least in low rank, see section (\ref{acm-semistable}).
In particular we classify these bundles up to rank $3$.

The manifold $V_5$ is a Fano threefold of {\em index two}, i.e. $\omega_{V_5} \simeq  \OO_{V_5}(-2)$,
where $\OO_{V_5}(1)$ is the minimal very ample line bundle on $V_5$. This is clear since
$V_5$ is obtained cutting $\G(\p^1,\p^4)$ with a generic $\p^6 \subset \p^9$.
\begin{equation} \label{G14-section}        
    V_5 = \G(\p^1,\p^4) \cap \p^6 \subset \p(\wedge^2 V)
\end{equation}

A quadric section of $V_5$ is a $K3$ surface, whose linear section is a {\em
canonical} curve of genus $6$, that is a curve embedded by the canonical
linear system. Such curve is {\em generic}. This is related
to the classification of Fano threefolds and to the moduli space of $K3$
surfaces, see \cite{mukai:curves-K3}.

From the point of view of birational geometry, $V_5$ is the
contraction of the proper preimage of a $2$-dimensional quadric in $Q_3$
under the blow up of $Q_3$ with center at a rational normal cubic. For a more detailed
description of many properties of this and other threefolds,
one may consult \cite{fano-encyclo}.

\section{Quasi--homogeneous structure} \label{quasi-homogeneous}

The linear orbit of a $d$-tuple of points in $\p^1=\p^1_{\kk}=\p(Y)$ is the orbit 
of a set of $d$ points in $\p^1$ under the standard action of $\SL(2)=\SL(Y)$.
The orbit sits in $\p(\Sym^d Y)$
where $\SL(2)$ naturally acts by the weight-$d$ representation, denoted
by $Y_d$. We have $Y_d \simeq \Sym^d Y$.

The orbit has dimension at most $3$, and there are few cases in which
its closure is smooth. In this case the closure is actually a Fano threefold
with $b_3 = 0$ and the possible cases are
        \begin{enumerate}
        \item $\p^3 \simeq \p(Y_3)$ orbit of $x^3 + y^3$;
        \item $Q_3 \subset \p(Y_4)$ orbit of $x^4 + xy^3$;
        \item $V_5 \simeq \G(\kk^2,Y_4) \cap \p(Y_6)$ orbit of $x^5y - xy^5$;
        \item $U_{22} \subset \p(Y_{12})$ orbit of $xy(x^{10} + 11x^5y^5-y^{10})$.
        \end{enumerate}

This is studied in detail in \cite{aluffi-faber:P1-orbits},
\cite{mukai-umemura}, \cite{mukai:fano-3-folds}. Clearly the stabilizer is a finite
group and looking at isometries of the roots of such polynomials one checks that
it is isomorphic respectively to $S_3$, $A_4$, $S_4$, $A_5$.

\begin{prop}[Mukai--Umemura]
The manifold $V_5$ is isomorphic to the closure of $\SL(2) \cdot (x^5y -xy^5)$ in
$\p(Y_6) = \p ^6$. 
The boundary divisor has an open part given by $\SL(2) \cdot
x^5y$ whose complementary is the degree-$6$ rational normal curve $\SL(2)\cdot
x^6$. 
\end{prop}

On the other hand, let $B$ be a $3$-dimensional vector space,
$\p^2=\p(B)$ and the dual plane $\pd^2=\p(B^*)$. Denote by $R = \Sym
B$ the coordinate ring, and let $S = \Sym^* B$ act by apolarity on $R$.
Write $z_0,z_1,z_2$ for the generators of $B$ and $\partial_0,\partial_1,\partial_2$
for those of $B^*$. Let $F\in \Sym^2 B$ be the equation of a smooth  conic in $\p(B)$
and $V = \Sym^2B / (F)$. 

\begin{prop}[Mukai]
Let $\VPS(3,\V(F))$ be the variety of polar $3$-sides
to the conic $\V(F)$, i.e. the closure in $\Hilb_3(\pd^2)$ of the set
        $$\{(f_0,f_1,f_2) \in \Hilb_3(\pd^2) | f_0^2 + f_1^2 + f_2 ^2= F\}$$

Then the variety $V_5$ is isomorphic to $\VPS(3,\V(F))$.
\end{prop}

The conic $F$ endows $B$ with a natural $\SO(3)$-action,
and we consider the dual conic $F^{-1}$ as an element of $\Sym^2B^*$. 
Taking $2 \times 2$ minors of the natural evaluation $B^* \ts \Sym^2B \rt B$ provides a map 
$\wedge^2 B^* \simeq B \sr \wedge ^2 V^*$ and the 
linear section (\ref{G14-section}) in $\p(\wedge^2 V)$ corresponds to the net of
alternating forms $\sigma \colon \wedge^2 B^* \simeq B \sr \wedge ^2 V^*$ described by the following
        $$\xymatrix@R-4ex@C-3ex{
        {\wedge^2 B^* \ts \wedge^2 V} \ar[r] & {\Sym^2 B} \ar[r]^-{F^{-1}} & {\kk} \\
        {\partial_1 \wedge \partial_2 \ts f_0 \wedge f_1} \ar@{|->}[rr] && {F^{-1}}(\partial_1 \wedge \partial_2 ( f_0 \wedge f_1))
        }$$

The standard $2:1$ covering $\phi \colon \SL(2) \rt \SO(3)$ gives
$B$ a structure of weight-$2$ representation, and we get $B \simeq Y_2$ and $V\simeq Y_4$.
Moreover the arrow $\sigma \colon B \rt \wedge ^2 V^*$ giving the net of
alternating forms is equivariant, i.e. is just the kernel of the projection 
$\wedge^2V^* \rt \Sym ^6 Y \simeq Y_6$. So $\p^6\simeq\p(Y_6)$.

If we choose $F$ to be given by the identity matrix, then
$(z_0^2 ,z_1^2, z_2^2)/(F)$ represents a polar triangle of
$C=\V(F)$. Its $\SO(3)$-orbit is $\SO(3) \cdot z_0^2 \wedge z_1^2$
in $\wedge ^2 V$ and under a $\phi$-equivariant projection $z_0^2 \wedge
z_1^2$ corresponds to $(x^5y -xy^5)$.

Now let us define the variety $\G(2\times 3;B^*)$ given by 
classes of $2\times 3$ matrices over the space $B^*$, that is
    $$\G(2\times 3;B^*) = \p(\M_{2,3}(B^*))//\SL(2)\times \SL(3)$$

The variety $\G(2\times 3;B^*)$ compactifies the open subset  $\Hilb_3(\pd^2)^{\circ}$ of 
triples of distinct points and is endowed with the natural bundles $\tilde{Q}^*$ 
and $\tilde{U}$ respectively of rank $3$ and $2$. This variety has been deeply 
studied by Drezet in \cite{drezet:hauteur-nulle}, where he proves
(\cite[Theorem 4]{drezet:hauteur-nulle}) that such compactification
actually provides a morphism
      $$\Hilb_3(\pd^2) \rt \G(2\times 3, B^*)$$
which is the blow up of $\G(2\times 3, B^*)$ at $\p^2 \subset \G(2\times 3, B^*)$. 
The bundles $\tilde{Q}^*$ and $\tilde{U}$ extend the sheaves over 
$\Hilb_3(\pd^2)^{\circ}$ whose fiber over $Z$ is given
respectively by the ideal generators and the first order syzygies
of the subscheme $Z \subset \pd^2$. It follows that $\HH^0(\tilde{Q})^* \simeq \Sym^2 B^*$
and $\HH^0(\tilde{U})^* \simeq \ker (m : \Sym^2 B^* \ts B^* \sr \Sym^3 B^*)$, where
$m$ is the multiplication in $S$, the coordinate ring of $\pd^2$.
The restrictions of $\tilde{U}$ and $\tilde{Q}$ to $V_5$ coincide with $U$ and $Q$, defined
as restrictions from $\G(\kk^2,V)$.

\begin{lem}
Let $F\in \Sym^2 B$ be a smooth conic and consider $F$ as a function
$F \colon \Sym^2B^* \sr k$. Then $V_5$ is isomorphic to
the following subset of $\Hilb_{3}(\pd ^2)$
       $$\{Z \in \Hilb_{3}(\pd ^2) | \HH^0(I_Z(2)) \subset \ker F \}$$

Moreover $V_5$ is isomorphic to 
the zero locus of the section $F\in \HH^0(\G(2\times 3;B^*),\tilde{Q}))$
\end{lem}
\begin{proof}
One may check this when $F = z_0^2+z_1^2+z_2^2$, taking
the polar triangle $Z\in \Hilb_{3}(\pd ^2)$ given by the three lines $(\V(z_0),\V(z_1),\V(z_2))$
and then using the action of $\SO(3)$, indeed the section $F$ is invariant under $\SO(3)$. 
The ideal $I_Z$ of such a point in $S = \Sym B^*$ is generated by $(\partial_0 \partial_1,
\partial_0 \partial_2, \partial_1 \partial_2)$. Clearly such
generators lie in $\ker F$. Notice also that the section $F$ vanishes in the
expected codimension. Furthermore the total spaces of $Q$ and $U^*$ in this framework are
    \begin{align*}
    \HH^0(V_5,Q)^*   \simeq & \,\, V^* \simeq \ker (F \colon \Sym^2 B^* \sr \kk) \\
    \HH^0(V_5,U^*)^* \simeq & \,\, V \simeq \ker (m \colon V^* \ts B^* \sr \Sym^3 B^*)
      \end{align*}

We can identify the kernel
of $m$ with $V\simeq(V^*)^*$ because, denoting by $J_F$ the ideal generated by $V^*$, 
$S/J_F$ is a codimension three Gorenstein ring (the {\em apolar} ring to $F$),
whose resolution is given by the structure theorem in \cite{buchsbaum-eisenbud:algebra-structures}.
Such resolution provides
    $$V^* \simeq \Tor_1^S(R/J_F,k)_2 \simeq \Tor_2^S(R/J_F,k)_3^* \simeq (\ker m)^*$$  

Notice also that all the arrows in this setup are $\SL(2)$-invariant.
\end{proof}

\section{Vector bundles and mutations} \label{bundles}

Here we consider bundles over $V_5$. An exceptional collection of
bundles over $V_5$ is known after \cite{orlov:V5}, where helix conditions are also
proved. We implicitly refer to \cite{bondal:helices},
\cite{nogin:helices-on-fano} for implications of this. We will resketch
Orlov's exceptional collection here.

\vspace{0.3 cm}

Given two bundles $E$ and $F$
define $i_{E,F}$ as the dual evaluation $i_{E,F} : E \sr \Hom(E,F)^*\ts F$ and $\ev_{E,F}$ the
evaluation $\ev_{E,F}: \Hom(E,F)\ts E \sr F$. Moreover, denote by $c_i \in \Z$ the
$i$-th Chern class of a sheaf $E$, meaning $c_i(E)=c_i \xi_i$, where $\xi_1 = c_1(\OO_{V_5}(1))$ and
$\xi_2$ (resp. $\xi_3$) is the class of a line (resp. of a point) in $V_5$.

\begin{dfn}
Given two bundles $E$ and $F$, if $i_{E,F}$ is injective we define the right mutation
$\Rm_F(E) = \cok (i_{E,F})$. If $\ev_{E,F}$ is surjective we define the left mutation $\Lm_{E}(F) = \ker (\ev_{E,F})$.
We refer the reader to \cite{gorodentsev:mutations},
\cite{rudakov:helices}, \cite{drezet:beilinson} for more general definitions
and essential properties of mutations.
\end{dfn}

We consider the following vector bundles over $V_5$. The bundle $U$ is the 
universal subbundle with $\rk(U)=2$, $c_1(U)=-1$ and $c_2(U)=2$.
The bundle $Q$ is the universal quotient bundle, $\rk(Q)=3$, $c_1(Q)=1$, $c_2(Q)=3$, $c_3(Q)=1$.

\begin{dfn} \label{collection}
Define the collection $\big( G_3,\ldots, G_0 \big)$ as $\big(  \OO(-1),U,Q^*,\OO \big)$.
Define also the dual collection
$\big( G^3,\ldots, G^0 \big)$ as $\big(\OO(-1),\wedge^2 Q^*,U,\OO \big)$.
\end{dfn}

The following
lemma is straightforward.
\begin{lem} \label{basic-bundles}
We have the exact sequence with $\SL(2)$-equivariant maps
        \begin{equation} \label{universal}
        0 \rt U \xrightarrow{i_{U,\OO}} V \ts \OO \rt Q \rt 0
        \end{equation}
and we have the following isomorphisms of $\SL(2)$-modules
\begin{equation} \label{V5-sections}
        \HH^0(\OO(1)) \simeq Y_6 \qquad \HH^0(U^*) \simeq V^* \simeq Y_4 \qquad \Hom(U,Q^*) \simeq B \simeq Y_2
        \end{equation}

Moreover we have the descriptions of the Hilbert scheme of lines and
conics in $V_5$
        \begin{align}
        & \Hilb_{2t +1} (V_5) \simeq \p(\HH^0(U^*)) \simeq \p(Y_4) = \p^4 \\
        & \label{castelnuovo-nuovo} \Hilb_{t +1} (V_5) \simeq \p(B) = \p^2
        \end{align}
\end{lem}

The vector spaces in (\ref{V5-sections}) inherit an invariant duality,
hence we will deliberately confuse them with their dual from now on.
Notice that the net $\sigma$ of alternating forms on the space $V$ is
just a $3$-subspace of $\HH^0(\wedge^2 U_{\G}^*) \simeq
\HH^0(\OO_{\G}(1)) \simeq \wedge^2 V^*$.  On other hand, given an
element $b$ in $B\simeq \Hom(U,Q^*)$, the homomorphism $b : U \sr Q^*$ takes $u \in U$
to $\sigma_b(u,-) \colon Q \sr \kk$ with $\sigma$ described above.  It
is elementary to check the degeneration locus of this homomorphism and 
write the exact sequence
  \begin{equation} \label{ideal-of-line}
  0 \rt U \overset{b}{\rt} Q^* \rt I_L \rt 0
  \end{equation}
where $L$ is a line in $V_5$. We thus recover (\ref{castelnuovo-nuovo}), 
actually known to Castelnuovo, \cite{castelnuovo:retta-quattro-dimensioni}.

\begin{prop} \label{basic-mutations}
  There are $\SL(2)$-equivariant exact sequences
        \begin{align} 
        \label{QQ}  & 0 \rt Q(-1) \xrightarrow{i_{Q(-1),U}} Y_2 \ts U \xrightarrow{\ev_{U,Q^*} }Q^* \rt 0 \\
        \label{E_9} & 0 \rt \OO(-1) \xrightarrow{i_{\OO(-1),U}} Y_4 \ts U \rt E_9 \rt 0  \\
        \label{omega-helix}    & 0 \rt E_9 \rt Y_4 \ts Q^* \rt \Omega_{\p^6}^1(1)_{|V_5} \rt 0
        \end{align}
        where $E_9$ is a rank-9 vector bundle. The sequences fit
        together to the helix below
        $$\xymatrix@-3ex{ {\OO(-1)} \ar@{^{(}->}[rr] && Y_4 \ts U
          \ar[dr] \ar[rr]
          && Y_4 \ts Q^*\ar[dr] \ar[rr] &&Y_6 \ts \OO \ar@{->>}[rr] && {\OO(1)}\\
          &&& E_9 \ar[ur] && {\Omega_{\p^6}^1(1)_{|V_5}} \ar[ur] }$$
\end{prop}

\begin{proof}
  First recall that $\wedge^2 Q^* \simeq Q(-1)$ and $B\simeq
  \Hom(U,Q^*) \simeq Y_2$. Clearly the evaluation $Y_2 \ts U \rt Q^*$
  is $\SL(2)$-equivariant and $\wedge^2 Q^* \simeq Q(-1)$ so by self
  duality we have the exact sequence with invariant
  maps (\ref{QQ}) where the first map in coincides with $i_{Q(-1),U}$.

  First notice that $\Ext^1(Q,U(1))=0$. Indeed one can compute  
$\HH^1(V_5,Q^*\ts U(1))=0$ taking global sections in the
Koszul complex of $V_5 \subset \G(k^2,V)$ and using
$\HH^p(\G(k^2,V),Q^*\ts U(2-p)) = 0$ for $1 \leq p \leq 4$.

Then applying $\Hom(-,U(1))$ to the sequence (\ref{universal}),
and using the identifications $\Hom(U,U(1))\simeq \Hom(U(-1),U)$
and $\Hom(Q,U(1)) \simeq \Hom(Q,U^*) \simeq \Hom(U,Q^*) \simeq B$
we get
        \begin{equation} \label{Hom(U(-1),U)}
        0 \rt B \xrightarrow{\,\, \sigma \,\,} V^* \ts V^* \rt \Hom(U(-1),U) \rt 0
        \end{equation}

Since there is a unique (up to scalar) equivariant map $B \sr V^* \ts V^*$, it
coincides with $B \xrightarrow{\, \sigma \,} \wedge^2 V^* \sr V^* \ts V^*$, so we denote it by $\sigma$
by abuse of notation. Since $Y_4\ts Y_4 \simeq Y_8 \oplus Y_6 \oplus Y_4 \oplus Y_2 \oplus Y_0$,
the group $\Hom(U(-1),U)$ is a rank-$22$ vector space isomorphic to $Y_8
\oplus Y_6 \oplus Y_4  \oplus Y_0$. 
We get a rank-$42$ bundle by
the following right mutation (i.e. $E_{42}=\Rm_{U}(U(-1))$)
        \begin{equation} \label{E_42}
        0 \rt U(-1) \rt \Hom(U(-1),U) \ts U \rt E_{42} \rt 0
        \end{equation}

Mutating again gives the following the exact sequence (\ref{E_9}). One
computes $\Hom(E_9,Q^*) \simeq V \simeq Y_4$ and sees that $E_{42}$ is also the
left mutation of $Q^*$ with respect to $E_9$, i.e. we have an exact sequence
        \begin{equation} \label{E_42E_9}
        0 \rt E_{42} \rt Y_4 \ts E_9 \rt Q^* \rt 0
        \end{equation}

Finally denote by $E_6$ the bundle given by right mutating $E_9$ with
respect to $Q^*$
        $$0 \rt E_9 \rt Y_4 \ts Q^* \rt E_6 \rt 0$$

Since $\Omega_{\p^6}^1(1)_{|V_5}$ is obtained as the kernel
        $$0\rt \Omega_{\p^6}^1(1)_{|V_5} \rt \Hom(\OO,\OO(1)) \ts \OO \simeq Y_6 \ts \OO \rt \OO(1) \rt 0$$
we get $E_6 \cong \Omega_{\p^6}^1(1)_{|V_5}$ thus obtaining (\ref{omega-helix}). The
cycle of left mutations in the exceptional collection $\big(U,Q^*,\OO,\OO(1) \big)$
then gives the helix condition for $\OO(1)$, proving the proposition. This appears
also in Orlov's paper \cite{orlov:V5}.
\end{proof}

\section{Resolution of the diagonal} \label{resolution}

We will write a resolution of the diagonal over $V_5$ in terms of
the collection of Definition \ref{collection}. We will make use of
the terminology of derived categories and derived functors, and we refer to
\cite{gelfand-manin:homological}
for definition and basic properties of these categories.
So let $\D^b(V_5)$
be the bounded derived category of coherent sheaves on $V_5$. An object
in $\D^b(V_5)$ is represented by a complex $\KK$ with finite nonzero cohomology.
We write $i$-th term of $\KK$ as $\KK^i$ and
we can shift $\KK$ setting $\KK[p]^j=\KK^{j+p}$. Also we set $\KK^* = \Rr \mathcal{H}om(\KK,\OO)$.

\begin{thm} \label{V5-main} The variety $V_5$ admits the following resolution of the diagonal:
        $$ 0 \rt \OO(-1,-1) \rt U \boxtimes \wedge^2 Q^*
                \rt Q^* \boxtimes U \rt \OO \rt \OO_{\Delta} \rt 0 $$
where the arrows are determined as the unique (up to scalar) morphisms invariant for the $\SL(2)$-action.
Denote by $d_i : G_{i+1} \bx G^{i+1} \sr G_i \bx G^i$ the $-i$-th
differential $\CC_{\Delta}^{-i-1} \sr \CC_{\Delta}^{-i}$ in the complex $\CC_{\Delta}$ over the diagonal.
Then we have
        $$\CC_{\Delta}^*[3] \simeq  \tau^*(\CC_{\Delta}(-1))$$
where $\tau$ is the involution that interchanges factors in $V_5\times
V_5$.
\end{thm}
\begin{proof}
First notice that such morphisms are indeed unique up to scalar, since they are nonzero and 
\begin{align*}
  & \Hom(\OO(-1,-1),U \boxtimes \wedge^2 Q^*) \simeq V^* \otimes V \simeq \End(Y_4)\\
  & \Hom(U \boxtimes \wedge^2 Q^*,Q^* \boxtimes U) \simeq B^* \otimes B \simeq \End(Y_2)\\
  & \Hom(Q^* \boxtimes U,\OO) \simeq V \otimes V^* \simeq \End(Y_4)
\end{align*}

In each case the differential we have defined is the identity element respectively over 
$V$ or $B$ or $V$, which is the unique $\SL(2)$-invariant element in the endomorphism groups,
since the tensor product decomposition of irreducible $\SL(2)$-modules contains a unique trivial
summand. On the other hand 
$$
  \Hom(\OO(-1,-1),Q^* \boxtimes U) \simeq Y_6 \otimes Y_4 \qquad
  \Hom(U \boxtimes \wedge^2 Q^*,\OO) \simeq Y_4 \otimes Y_6
$$

These modules contain no nontrivial invariant element, hence the composition of invariant maps must be zero.

To prove that the complex is actually exact, we need a more explicit
description of such maps. We will provide a fiberwise description
over a point $(p,q) \in V_5 \times V_5$.  However, the duality of the
complex tells us that we need only to prove the exactness in $\OO$ and
$U \boxtimes \wedge^2 Q^*$. Recall from Lemma \ref{basic-bundles} and Proposition \ref{basic-mutations}
the maps $i_{U,\OO}$ and $i_{\OO(-1),U}$ and consider $i_{Q(-1),U}$ $i_{Q^*,\OO}$ and
$i_{\OO(-1),Q(-1)}$. Then we have the following commutative diagrams
        $$\xymatrix@C+2ex{d_0 : Q^* \boxtimes U \rt \OO &
          Q^*_q \ts U_p\ar[dr]^{d_0} \ar_-{i_{Q^*,\OO} \boxtimes i_{U,\OO}}[d]       &\\
                                &       V \ts V^*\ar[r]_-{\chi}         & {\kk}}$$

        $$\xymatrix@C+2ex{ d_1 :U \boxtimes \wedge^2 Q^*
                \rt Q^* \boxtimes U &
        U_p \ts \wq_q \ar[dr]^{d_1} \ar[d]_-{1 \boxtimes i_{Q(-1),U}}      &\\
        & U_p \ts B \ts U_q \ar[r]_-{\ev_{U,Q^*} \ts 1}           & Q^*_q \ts U_p}$$

        $$\xymatrix@C+3ex{ d_2 : \OO(-1,-1) \rt U \boxtimes \wedge^2
        Q^* & 
        {\OO_p(-1) \ts \OO_q(-1)} \ar[dr]^{d_2} \ar[d]_-{i_{\OO(-1),U}\boxtimes i_{\OO(-1),Q(-1)}}            &\\
        & U_p \ts V \ts V^* \ts \wq_q \ar[r]_-{1\ts \chi \ts 1}   & U_p \ts \wq_q}$$

where the morphism $\chi$ here is the natural evaluation $V^* \ts V \sr k$.
The definition of $d_0$ immediately tells that the complex
is exact in $\OO$ since the evaluation vanishes only when $U_p =
Q^*_q$ i.e. for $p=q$. To prove the other step it is useful to give another description of $d_2$
        $$\xymatrix@C+4ex{
        {\OO_p(-1) \ts \OO_q(-1)} \ar@{-->}[r]_-{d_2} \ar[d]_-{i_{\OO(-1),U} \boxtimes i_{\OO(-1),U}}
          & U_p \ts \wq_q \ar[d]^-{1 \boxtimes i_{Q(-1),U}}\\
        U_p \ts V \ts V \ts U_q \ar[r]_-{1\ts \sigma^* \ts 1}       & U_p \ts B \ts U_q}$$
where $\sigma ^* \colon V \ts V \rt B \simeq B^*$ is the transpose of the map
$\sigma$ defined in section $(\ref{quasi-homogeneous})$, composed with
$V \ts V \rt \wedge^2 V$.

All such definitions coincide up to scalar with the previous ones by uniqueness of the
invariant map since they produce nonzero maps. From these descriptions one readily gets
        $$\ker (d_1) \simeq \wedge^2 Q^*_p \ts U_q \cap U_p \ts \wedge^2 Q^*_q \subset U_p \ts B\ts U_q$$

On the other hand the image of $d_2$ is the image in $U_p \ts B \ts
U_q$ of $1\ts \sigma \ts 1$ restricted to $\OO_p(-1) \ts \OO_q(-1)$. This coincides with $\ker(d_1)$ if
\begin{align*}
        & \wedge^2 Q^*_p \ts U_q = \im (1\ts \sigma^* \ts 1) && \mbox{restricted to
        $\OO_p(-1) \ts V \ts U_q$} \\
        & U_p \ts \wedge^2 Q^*_q = \im (1\ts \sigma^* \ts 1) && \mbox{restricted to
        $U_p \ts V \ts \OO_q(-1)$}
\end{align*}

Consider the first equality (they are clearly symmetric) and notice that we may factor out the
identity over $U_q$ and thus reduce it to
        $$ \wedge^2 Q^* = \im (1\ts \sigma^*) \qquad \mbox{restricted to
        $\OO(-1) \ts V$}$$

This is provided by the commutativity of the following diagram, since the only equivariant
morphism $\OO(-1) \ts V \sr \wedge^2 Q^*$ is the map given by the universal quotient (twisted by $-1$), which is of
course surjective.
        \begin{equation} \label{key}
        \xymatrix@R-1ex{
                        &               & 0\ar[d]                                       & 0\ar[d]                               & \\
        0 \ar[r]        &U(-1) \ar[r]^-{i_{U,\OO}}   & V \ts \OO(-1) \ar[r]^-{\ev_{\OO,Q}}\ar[d]_-{1 \ts i_{\OO(-1),U}}     & {\wedge^2 Q^*} \ar[d]^{i_{Q(-1),U}}  \ar[r] & 0\\
                        &               & V \ts V \ts U \ar[r]_-{\sigma^*\ts 1}         & B \ts U \ar[d]^{\ev_{U,Q^*}}  \ar[r]          & 0\\
                        &               &                                               & Q^*   \ar[d]                          &\\
                        &               &                                               & 0                                     &}
        \end{equation}
\end{proof}

The immediate consequence of this theorem can be formulated at the level of derived categories.
Indeed the functor ${\Rr
p_2}_*(\CC_{\Delta} \overset{\Ll}{\ts} p_1^*(-))$ associated to the complex
$\CC_{\Delta}$ resolving the diagonal (where $p_1,p_2$ are the two
projections $V_5\times V_5 \sr V_5$) is isomorphic the identity in
$\D^b(V_5)$. Then for any sheaf $F$ we
have a complex $\CC_F$, exact except in cohomological degree $0$, whose
cohomology is $F$. The $i$-th term of $\CC_F$ is the sum $\oplus_j \HH^{j+i}(F \ts G^j) \ts G_j$.
In order to construct the $i$-th term of $\CC_F$ it is thus sufficient to consider
$i$-th diagonal (where the $0$-th diagonal is the middle one and the $3$rd (resp. $-3$rd)
is the upper--right (resp. lower--left) corner)
of the following cohomology matrix
\begin{equation} \label{beilinson-table}
\left(
\begin{aligned}
        \hh^3 (F \ts \OO(-1))   \quad && \hh^3 (F \ts Q(-1))    \quad && \hh^3 (F \ts U)        \quad && \hh^3 (F)  \\
        \hh^2 (F \ts \OO(-1))   \quad && \hh^2 (F \ts Q(-1))    \quad && \hh^2 (F \ts U)        \quad && \hh^2 (F)  \\
        \hh^1 (F \ts \OO(-1))   \quad && \hh^1 (F \ts Q(-1))    \quad && \hh^1 (F \ts U)        \quad && \hh^1 (F)  \\
        \hh^0 (F \ts \OO(-1))   \quad && \hh^0 (F \ts Q(-1))    \quad && \hh^0 (F \ts U)        \quad && \hh^0 (F) \\
\end{aligned}
\right)
\end{equation}
We say that the collection $(G_3,\ldots,G_0)$ of Definition \ref{collection} is a {\em basis} for $\D^b(V_5)$
and that $(G^3,\ldots,G^0)$ is its {\em dual basis}.
Orlov in \cite{orlov:V5} proves that the exceptional sequence $(U(-1),\OO(-1),U,\OO)$
generates $\D^b(V_5)$ since it satisfies helix conditions.

\begin{corol} \label{beilinson-V5}
  Any coherent sheaf $F$ on $V_5$ is isomorphic to the cohomology of a
  complex $\CC_F^{\bullet}$ (dually, of a complex $\DD_F^{\bullet}$)
  whose terms are given by
        $$\CC_F^k = \bigoplus _{i-j = k} \HH^i(F \ts G^j) \ts G_j \qquad
        \DD_F^k = \bigoplus _{i-j = k} \HH^i(F \ts G_j) \ts G^j $$

The duality in Theorem \ref{V5-main} implies $\DD_F^{\bullet} \simeq \CC_F^{\bullet}(-1)^*$.
\end{corol}

\begin{rmk}
  The construction of $V_5$ as subvariety of $\G(2\times 3;B^*)$
  allows to prove Theorem \ref{V5-main} without making use of the
  $\SL(2)$-action.  The diagram (\ref{key}) can be completed in the
  following, that summarizes the various mutations (\ref{E_42}),
  (\ref{E_9}), (\ref{E_42E_9}).
  $$
  \xymatrix@C-2ex{
    U(-1) \ar[r] \ar[d]                     & V \ts \OO(-1) \ar[r]\ar[d]            & {\wedge^2 Q^*} \ar[d] \\
    \Hom(U(-1),U) \ts U \ar[r]\ar[d]        & V \ts V \ts U \ar[r] \ar[d]           & B \ts U \ar[d]        \\
    E_{42} \ar[r] & V \ts E_9 \ar[r] & Q^* }$$

  Rows and columns are
  exact and zeroes all around the diagram are omitted for brevity.
  The arrows in the central row are just the the sequence
  (\ref{Hom(U(-1),U)}) tensorized with $1_U$ and all maps are invariant
  evaluations.
\end{rmk}

\section{Splitting criterion} \label{splitting-section}

We have the following well--known splitting criterion (cfr. \cite{ottaviani:horrocks}).
\begin{prop} \label{splitting-general} Let $X$ be a connected projective variety embedded by $\OO_X(1)$
and containing a straight line $L$ with ideal sheaf $I_L$. Let $E$ be a vector bundle on $X$. Then
$E$ splits as a sum of line bundles if:
        \begin{equation} \label{split2}
        \HH^1(E\ts I_L \ts \OO_X(1)^{\ts t}) = 0 \quad \mbox{for any t $\in \Z$} 
        \end{equation}
\end{prop}
\begin{proof} 
  Restricting $E$ to $L$ gives
  $$
  0 \rt E \ts I_L \rt E \rt E_{|L} \rt 0$$

  Let $F$ be the split
  bundle on $X$ such that $F_{|L} \simeq E_{|L}$. Then consider the
  exact sequence
  $$
  0 \rt F^* \ts E \ts I_L \rt F^* \ts E \rt F_{|L}^* \ts E_{|L} \rt
  0$$

  Condition (\ref{split2}) gives $\HH^1(F^* \ts E \ts I_L) = 0$,
  hence surjectivity in the map
  $$
  \HH^0(F^* \ts E) \rt \HH^0(F_{|L}^* \ts E_{|L}) $$
  so the
  identity of $F_{|L}$ lifts to a morphism of $F$ to $E$. The
  determinant of such a morphism is a constant\ since it lies in
  $\HH^0(\OO)$. This constant is nonzero by restriction to $L$, i.e.
  $E$ is indeed isomorphic to $F$.
\end{proof}

\begin{prop} \label{splitting} A vector bundle $E$ on $V_5$ splits if and only if, for
any $t\in \Z$ 
\begin{enumerate}[i)]
\item $\HH^1(E\ts Q^*(t)) = 0$
\item $\HH^2(E\ts U(t)) = 0$
\end{enumerate}
\end{prop}

\begin{proof} A sum of line bundles obviously satisfies the
conditions since $\hh^2(U(t)) = 0 = \hh^1(Q^*(t))$ for any
$t$.

Viceversa, consider the ideal sheaf of a line in $V_5$ as the cokernel of
a map $b : U \sr Q^*$ as in the exact sequence (\ref{ideal-of-line}).
Then the hypothesis implies the condition of Proposition \ref{splitting-general}.
\end{proof}

This proposition can be thought of as an analog of the Horrocks splitting
criterion on the projective space $\p^n$ (see \cite{horrocks:punctured}), namely that a
bundle splits if and only if its cohomology modules other than $0$ and
$n$ vanish in all degrees.
Next we will show how to deduce a generalization of Proposition \ref{splitting} from the Corollary
\ref{beilinson-V5}. By applying convenient mutations to the exceptional collection
$\big(  \OO(-1),U,Q^*,\OO \big)$, we will find another resolution of the diagonal
over $V_5 \times V_5$, suitable to deduce the splitting criterion.

\begin{lem} \label{mutated-V5}
The exceptional collection $\big( G_3,\ldots, G_0 \big) = \big(  \OO(-1),U,Q^*,\OO \big)$ can be
mutated to the exceptional collection
        $$\big(  U(-1),Q^*(-1),\OO(-1),\OO \big)$$
\end{lem}
\begin{proof}
This is elementary. The first step is
        $$\big(  \OO(-1),U,Q^*,\OO \big) \mapsto \big(  \OO(-1),\wedge^2 Q^*,U,\OO \big)$$
by mutating $U$ and $Q^*$. This is nothing but passing to the dual
collection. Next we mutate $\OO(-1)$ and $\wedge^2Q^*$ to get
        $$\big(  \OO(-1),\wedge^2 Q^*,U,\OO \big) \mapsto \big(  U(-1),\OO(-1),U,\OO \big)$$
finally we mutate $\OO(-1)$ and $U$ and get
        $$\big(  U(-1),\OO(-1),U,\OO \big)\mapsto \big(  U(-1),Q^*(-1),\OO(-1),\OO \big)$$
\end{proof}

\begin{corol} \label{Psi-V5}
The variety $V_5$ admits the mutated resolution of the diagonal
        $$0 \sr U(-1)\bx U \sr Q^*(-1) \bx Q^* \sr \OO(-1) \bx \Omega_{\p^6}^1(1)_{|V_5} \sr \OO \sr \OO_{\Delta} \sr 0$$
\end{corol}

\begin{proof}
It suffices to perform the left mutations of the collection $\big(
G_3,\ldots, G_0 \big)$ as in Lemma \ref{mutated-V5} and
the appropriate right mutations in the dual collection $\big( G^3,\ldots, G^0 \big)$.
Notice that the first step in Lemma \ref{mutated-V5} allows to switch the basis $\big(
G_3,\ldots, G_0 \big)$ with the dual basis. We need the following right mutations
        $$\big(  \OO(-1),U,Q^*,\OO \big) \mapsto \big(  U,E_9,Q^*,U,\OO \big)$$
and
        $$\big(  U,E_9,Q^*,U,\OO \big) \mapsto \big(U,Q^*,E_6,\OO \big)$$

But we have seen in section (\ref{bundles}) that $E_6 \simeq
\Omega_{\p^6}^1(1)_{|V_5}$.
\end{proof}

This allows to give an algebraic proof and a generalization of the
splitting criterion \ref{splitting}, in the spirit of \cite{ancona-ottaviani:beilinson-quadrics}.

\begin{corol} \label{splitting-beilinson}
Let $E$ be a sheaf on $V_5$ such that 
$$\hh^1(Q^* \ts E(t)) = \hh^2(U \ts E(t)) = 0 \qquad \forall \, t\in \Z$$ 

Then $E$ is isomorphic to a sum of line bundles plus a sheaf $E'$
with $\dim(\supp(E'))=0$
\end{corol}

\begin{proof}
By Corollary \ref{Psi-V5}, the mutated cohomology table of a given bundle $E$
reads
\begin{equation*}
\left(
\begin{aligned}
        \star \qquad            && \star \qquad                 &&        \star \qquad \qquad           && \qquad \star \qquad \\
        \HH^2(U\ts E(-1))       && \star \qquad                 &&        \star \qquad \qquad           && \qquad \star \qquad\\
        \star \qquad            && \HH^1(Q^*\ts E(-1))          &&        \star \qquad \qquad           && \qquad \star \qquad\\
        \star \qquad            && \star \qquad                 &&        \qquad \HH^0(E(-1)) \qquad    && \qquad \HH^0 (E) \\
\end{aligned}
\right)
\end{equation*}

Suppose that the Hilbert polynomial of $E$ is not constant. Then there exists a $t_0$ such that $\hh^0(E(t_0)) \neq
0$ and $\hh^0(E(t_0-1)) = 0$. This and the conditions of the hypothesis
imply that all terms in the $-1$ diagonal are $0$
for $E(t_0)$. Then $\HH^0(E(t_0)) \ts \OO$ is a direct summand of
$E(t_0)$.

By induction on the leading term of the Hilbert polynomial of $E$, 
$E$ is a sum of line
bundles plus a sheaf $E'$ whose Hilbert polynomial is constant. Then $E'$ 
is supported at a finite set of points.
\end{proof}

By mutating again we get the collection $\big( \wedge^2
Q^*(-1) , U(-1),\OO(-1),\OO \big)$. However this gives the same splitting
principle since $\hh^1(U \ts E) = \hh^2(U(-1) \ts E^*)$ and $\hh^1(Q^*
\ts E) = \hh^2(\wedge^2 Q^*(-1) \ts E^*)$ and obviously $E$ splits if and only if
$E^*$ does. The following proposition is in the spirit of
\cite{arrondo-grana:G(1.4)}, namely one characterizes $U$ or $\OO$ by
certain cohomology vanishing.

\begin{prop} \label{UorO}
Let $E$ be an indecomposable sheaf over $V_5$ 
such that $\dim(\supp(E))>0$ and

\begin{enumerate}[i)]
        \item   $\HH^1(U \ts E(t)) = \HH^2(U\ts E(t)) = 0$ for all
        $t\in \Z$
        \item   $\HH^1(E(t)) = \HH^2(E(t)) = 0$ for all
        $t\in \Z$
\end{enumerate}

Then $E$ is
isomorphic to either $U(a)$ or to $\OO(a)$ for some $a \in \Z$.
\end{prop}

\begin{proof}
In case $E$ is a bundle 
one can prove this with a technique analogous to
\cite{arrondo-grana:G(1.4)}. However here we show a different
method. Let us consider the cohomology table of such $E$ for the {\em
dual basis} $\big(G^3,\ldots,G^0\big)$, i.e. we
consider the complex $\DD_E^{\bullet}$ of Corollary \ref{beilinson-V5}.

\begin{equation*}
\left(
\begin{aligned}
        \star \qquad            && \star \qquad         &&        \qquad \hh^2(Q^* \ts E) \qquad        && \qquad \star \qquad \\
        0 \qquad                && 0 \qquad             &&        \qquad \hh^2(Q^* \ts E) \qquad        && \qquad 0     \qquad\\
        0 \qquad                && 0 \qquad             &&        \qquad \hh^1(Q^* \ts E) \qquad        && \qquad 0 \qquad\\
        \star \qquad            && \star \qquad         &&        \qquad \hh^0(Q^* \ts E) \qquad        && \qquad \hh^0 (E) \\
\end{aligned}
\right)
\end{equation*}

Then, the bundle $\HH^1(Q^* \ts E)
\ts U$ is a direct summand of $\DD_E^0$ 
mapping to zero in $\DD_E^1$ and with zero differential from $\DD_E^{-1}$.
Therefore $E$ contains such bundle as a
direct summand. Since this happens for any twist of $E$, we conclude
that $\oplus_t \HH^1(Q^* \ts E(t)) \ts U(t)$
is a direct summand of $E$.

By factoring out such summand, we can suppose $\HH^1(Q^* \ts E(t)) =
0$ for all $t$. But this condition, together with $\HH^2(U \ts E(t))=0$ for
all $t$, implies that $E$ splits by Corollary \ref{splitting-beilinson}.

This proves one statement if $E$ is indecomposable. 
The converse follows taking cohomology of the symmetrized square of (\ref{universal})
        $$ 0 \rt \Sym^2 U \rt V \ts U \rt \wedge ^2 V \rt \wedge ^2 Q \rt 0$$
and using the facts
        $$\Sym^2 U \oplus \OO(-1) \simeq U \ts U \qquad \Sym^2 U^* \simeq \Sym^2U(2)$$

Then $\Sym^2 U$ is aCM as well al $U \ts U$. For the twist $U \ts U^*$
it just says that $U$ is exceptional.
\end{proof}

Similar considerations lead to the following two propositions.

\begin{prop}
An indecomposable sheaf $E$ over $V_5$ such that $\dim(\supp(E))>0$
is either $\OO(a)$
or $Q^*(a)$ if and only if it has no intermediate cohomology and
        $$\HH^1(Q \ts E(t) = \HH^2(Q \ts E(t)) = 0 \qquad \forall t\in \Z$$
\end{prop}

\begin{prop} \label{UorQorO}
Let $E$ be an indecomposable sheaf over $V_5$ such that $\dim(\supp(E))>0$ and
\begin{enumerate}[i)]
        \item   $\HH^1(U \ts E(t)) = 0$ for all
        $t\in \Z$
        \item   $\HH^1(E(t)) = \HH^2(E(t)) = 0$ for all
        $t\in \Z$
\end{enumerate}
Then $E$ is either $U(a)$ or $\OO(a)$ or $Q(a)$ for some $a \in \Z$.
\end{prop}

\begin{proof}
By looking at the cohomology table as in (\ref{UorO}), we remove the
direct summand $\oplus _t \HH^1(Q^* \ts E(t)) \ts U(t)$. 
Then using the dual of the universal sequence (\ref{universal})
tensorized by $E(t)$ we get $\hh^2(Q^* \ts E(t)) = 0$ for all $t$ by
the hypothesis.

Now looking again at the cohomology table we see that the module
$\oplus _t \HH^2(U \ts E(t)) \ts Q(t-1)$ is a direct summand of
$E$, since the only cokernel could sit in $\HH^2(Q^* \ts E) \ts U(t)$.
Again factoring out such summand we get that the remaining part
splits by Corollary \ref{splitting-beilinson}.
\end{proof}

\section{aCM semistable bundles} \label{acm-semistable}
In this section we focus on {\em aCM} bundles i.e. with the condition $\HH^p(V_5,E \ts \OO(t))=0$ for any $t \in \Z$ and $0<p<3 = \dim(V_5)$.
\begin{thm}
An indecomposable aCM sheaf $E$ over $V_5$ with $\dim(\supp(E))>0$ fits into the following exact sequence
        $$0 \rt E \rt F_U \oplus F_Q \oplus F_{\OO} \rt G_U\rt 0$$  
where the bundles $F_U$, $F_{Q}$, $F_{\OO}$ and $G_U$ are given by
        \begin{align*}
        & G_U = \bigoplus_{t\in \Z} \HH^1(U\ts E(t)) \ts U(1-t) \\
        & F_U = \bigoplus_{t\in \Z} \HH^1(Q^*\ts E(t)) \ts U(t) \\
        & F_Q = \bigoplus_{t\in \Z} \HH^2(U\ts E(t)) \ts Q(t-1) \\
        & F_{\OO} = \bigoplus_{t\in I} \HH^0(E(t)) \ts \OO(t)   \qquad
                I = \{t\in \Z|\hh^0(E(t))\neq    0 = \hh^0(E(t-1)) \}
        \end{align*}
\end{thm}

\begin{proof}
Since $\HH^1(U \ts E (t)) = \Ext^1(U(1-t),E)$ we have a universal extension
        $$0 \rt E \rt F \rt \oplus_{t} \HH^1(U \ts E (t)) \ts U(1-t) \rt 0$$

Now $F$ is still aCM because $E$ and $U$ are. Furthermore, since $U
\ts U$ is aCM and we have taken {\em all extensions} between $E$ and
$U$ we get
        $$\HH^1(F \ts U(t)) = 0 \qquad \mbox{for all $t \in \Z$}$$

Thus the result follows by (\ref{UorQorO}). To compute the expression
of $G_U$, $F_U$, $F_Q$ and $F_{\OO}$ it suffices to trace back the
proof of Proposition (\ref{UorQorO}).
\end{proof}

To sharpen the analysis, we now assume that $E$ is a {\em semistable bundle} (in the sense of
Mumford--Takemoto, with respect to the hyperplane divisor). Recall that for any
coherent sheaf $F$ one defines $\rk(F) = \rk(F^{\circ})$ where $F^{\circ}$ is
the restriction of $F$ to a Zariski open set where $F$ is locally free.
Recall also the definition of slope $\mu(F)=c_1(F)/\rk(F) \in \Q$, writing $c_1(F)$ as
an integer as in section (\ref{bundles}).

\begin{dfn}
A vector bundle $E$ on $V_5$ is semistable (resp. stable) with respect to $\OO(1)$ if, for any
coherent subsheaf $F \subset E$ with $0 < \rk(F) <\rk(E)$ we have $\mu(F) \leq \mu(E)$ (resp. we have $\mu(F) < \mu(E)$).
\end{dfn}
For details concerning semistable sheaves 
we refer for instance to \cite{huybrechts-lehn:moduli}. We first {\em normalize} that $E$, that is twist it until
        $$-r < c_1(E) \leq 0$$

Let us now read the cohomology table (\ref{beilinson-table}) of such
$E$. 
First notice that in the cohomology table of $E$ the last
column is zero. Indeed $\hh^0(E(-1)) = 0$ and
$\hh^3(E)=\hh^0(E^*(-2))=0$ by semistability. 
In the first column, the only nonzero term can be $\hh^0(E)$. 
Now recall that the tensor product of semistable bundles is semistable
(if they are both stable bundles one can use Hermite--Einstein metrics,
otherwise see \cite[Theorem 1.14]{maruyama:grauert}).
Thus $E \ts U$ and $E \ts Q$ are semistable, so we get
\begin{align*}
 & \hh^0(E \ts U) = 0 && \hh^3(E \ts U) = 0 \\
 & \hh^0(E \ts Q(-1)) = 0 && \hh^3(E \ts Q(-1)) = 0 \\
\end{align*}
\vspace{-0.8cm}

Then the cohomology table is
\begin{equation*}
\left(
\begin{aligned}
        0       && 0 \qquad \qquad      \qquad          && 0 \qquad                     &&  0 \quad \\
        0       && \qquad \hh^2(Q\ts E(-1))\qquad       && \hh^2(U \ts E)       &&  0 \quad \\
        0       && \qquad \hh^1(Q\ts E(-1))\qquad       && \hh^1(U \ts E)       &&  0 \quad \\
        0       && 0 \qquad \qquad      \qquad          && 0 \qquad                     &&  \hh^0 (E) \\
\end{aligned}
\right)
\end{equation*}

Then $E$ is the cohomology of the complex 
        $$\xymatrix@R-3ex@C-2ex{
                        &                                       & {\HH^2(Q \ts E(-1)) \ts U} \ar@{}[d]|{\oplus}\\
        0 \ar[r]        & {\HH^1(Q \ts E(-1)) \ts U} \ar[r]     & {\HH^1(U \ts E)} \ts Q^* \ar[r]       & {\HH^1(U \ts E) \ts Q^*} \ar[r]       & 0\\
                        &                                       & {\HH^0(U) \ts \OO} \ar@{}[u]|{\oplus}
        }$$

But in the above complex there are no nonzero maps bundles with the
same base i.e. there are no vertical arrows in the spectral sequence
associated to complex $\CC_{\Delta}$. This means that $E$ splits into the
direct sum of two bundles $E_1$ and $E_2$ defined by 
        \begin{align*}
                & 0 \sr {\HH^1(E \ts Q(-1))} \ts U \sr \HH^1(E\ts U) \ts Q^* \oplus \HH^0(E) \ts \OO \sr E_1 \sr 0\\
                & 0 \sr E_2 \sr \HH^2(E \ts Q(-1)) \ts U  \sr \HH^2(E \ts U) \ts Q^* \oplus \HH^0(E) \ts \OO \sr 0
        \end{align*}

If we look at the cohomology of $E$ a little more carefully, we see
that
        $$\hh^2(U \ts E) = \hh^1(U(-1) \ts E^*) = \hh^0(Q(-1) \ts E^*)$$
thus semistability tells that
        $$\hh^2(U \ts E) = 0 \qquad \mbox{if} \qquad c_1(E) >-\frac{2}{3} \rk(E)$$

Analogously we have $\hh^2(Q(-1) \ts E) = \hh^1(Q^*(-1) \ts E^*) = \hh^0(U \ts E^*)$ so that
        $$\hh^2(Q(-1) \ts E) = 0 \qquad \mbox{if} \qquad c_1(E) >-\frac{1}{2} \rk(E)$$
and notice that if a normalized bundle $E$ does not satisfy the above,
then $E^*(-1)$ does, except the case $c_1(E) = - \frac{1}{2}\rk(E)$.
On the other hand if $c_1(E) = - \frac{1}{2}\rk(E)$, the resolution of
$E$ provided above tells that $E$ is a sum of copies of $U$ plus the cokernel
of a map $U^a \sr (Q^*)^b$ for certain integers $a$ and $b$. But $2
c_1(E) = - \rk(E)$ implies $b = 0$, so that $E$ is a sum of copies of
$U$. Summing up we can state the following theorem

\begin{thm} \label{aCM-ss}
Up to dualizing and twisting, any aCM semistable
bundle $E$ of over $V_5$ is a sum of copies of $U$ plus the cokernel $E'$
        $$0 \rt U^a \rt (Q^*)^b \oplus \OO^c \rt E' \rt 0$$ 
and we have
\begin{enumerate}[i)]
\item $E$ stable and $\rk(E)>1$ imply $c=0$;
\item $c_1(E) \neq 0$ implies $c=0$;
\item $c_1(E) = - \frac{1}{2} \rk(E)$ implies $E \simeq \HH^1(Q^* \ts E) \ts U$;
\end{enumerate}
where the above integers are
        $$a = \hh^1(Q(-1)\ts E) \qquad b = \hh^1(U \ts E) \qquad c = \hh^0(E)$$
\end{thm}

\begin{corol}
The moduli space of semistable aCM bundles over $V_5$ with fixed invariants is
unirational. So is the versal deformation space of any aCM bundle over $V_5$.
\end{corol}

\begin{rmk} One would get the same result for a normalized bundle $E$
with nonzero degree replacing the semistability hypothesis with the assumption that
        $$\hh^0(E) = 0 \qquad \hh^3(E(-1)) = 0$$
\end{rmk}

\begin{rmk} \label{enrique-laura}
Arrondo and Costa classified rank-$2$ aCM bundles over Fano threefolds of index
$2$ in \cite{arrondo-costa} by Hartshorne--Serre correspondence. There are
three families of such bundles on any threefold as above, related
respectively to a line, a conic and a projectively normal elliptic
curve in the threefold. On $V_5$ we
have the following cases.
\begin{enumerate}[i)]
        \item \label{E-line} $E_L$, the strictly semistable bundle whose unique
        section vanishes along a line $L$. Then we clearly have
                $$M(2;0,1) \simeq \p^2$$
        \item $E_C$ the stable bundle given by a conic. It is exactly $U$ and has no moduli.
        \item \label{E-elliptic} $E_S$ the stable bundle corresponding to a degree-$7$
elliptic curve $S$. $c_1(E_S)=0$ and $c_2(E_S) =2$. The moduli space of $E_S$ is five--dimensional.
\end{enumerate}

The resolutions of these bundles are respectively
$$
        0 \rt U \rt Q^* \oplus \OO \rt E_L \rt 0 \qquad
0 \rt U^{\oplus 2} \rt Q^{* \, \, \oplus 2} \rt E_S \rt 0$$
\end{rmk}

The following lemma allows to recover Arrondo and Costa's classification of aCM
bundles in rank $2$ quoted in Remark (\ref{enrique-laura}).

\begin{lem}
An indecomposable aCM rank-$2$ bundle $E$ over $V_5$ is semistable. In particular,
the only families of such bundles are the ones listed in Remark (\ref{enrique-laura}).
\end{lem}

\begin{proof}
Take the first twist $t$ such that $E$ has sections. Consider such section
and its vanishing locus $Z$. The subvariety $Z$ cannot be empty because
$\hh^1(\OO(t))=0$ and $E$ is indecomposable. Moreover $Z$ cannot
have divisorial components since we took the first twist with
sections. Thus $Z$ is a curve in $V_5$ with ideal sheaf $I_Z$
        $$0\sr \OO \sr E(t) \sr I_Z \ts \OO(c_1(E(t))) \sr 0$$

Now since $E$ is aCM we find $\hh^1(I_Z)=0$. This implies that $Z$ is connected since
$\hh^0(\OO_Z) = 1+\hh^1(I_Z)$.
But if $E$ is not semistable the Chern class $c=c_1(E(t))$ must be strictly negative
so the dualizing sheaf of $Z$ is
        $$\omega_Z \simeq K_{V_5} \ts \det(E(t))_{|Z} \simeq \OO(-2+c)_{|Z}$$

Then the degree of the dualizing sheaf of the connected curve $Z$
would be strictly less than two, which is impossible.

The last claim follows from Theorem (\ref{aCM-ss}). Indeed if $c_1(E)=-1$ we get
$b=0$, so that $E \simeq U$. If $c_1(E)=0$ we have, besides $\OO \oplus \OO$, the two cases $c=1$,
$c=0$, corresponding respectively to (\ref{E-line}) and (\ref{E-elliptic}) of Remark (\ref{enrique-laura}).
\end{proof}

If we denote $A_S = \HH^1(U \ts E_S)$ and $B_S = \HH^1(Q\ts E_S(-1))$
then one has the exact sequence
        $$0 \rt B_S \rt B \ts A_S \rt \HH^1(Q^*\ts E_S) \rt 0$$
which gives the dual presentation of $E_S$
        $$0 \rt E_S \rt U^4 \rt Q(-1)^2 \rt 0$$
and we also recover the exact sequence
        $$0 \sr \kk \sr B_S \ts B_S \sr A_S \ts \HH^1(Q^* \ts E_S) \sr \HH^1(E_S \ts E_S) \sr 0$$
identifying the tangent to the moduli space. Notice that a degree-$7$
elliptic curve is given by an element of $\HH^0(E(1)) \simeq
\kk^{10}$. This agrees with the dimension (given by Riemann--Roch) of the Hilbert scheme
$\dim(\Hilb_{7t}(V_5))=14$, a $\p^9$ fibration over the $5$-dimensional
moduli space.
Since $\Ext^1(E_S,E_S) \simeq \kk^5$ there are nontrivial extensions
        $$0 \rt E_S \rt E_4^S \rt E_S \rt 0$$

Then we have a rank-$4$ aCM bundle $E_4^S$ and its resolution reads
        $$0 \rt U^{\oplus 4} \rt Q^{* \, \, \oplus 4} \rt E_4^S \rt 0$$

The map $U^{\oplus 4} \sr Q^{* \, \, \oplus 4}$ deforms the
diagonal map $(\varphi,\varphi)$ where $\varphi \colon U^{\oplus 2} \sr Q^{* \, \, \oplus 2}$. 

Finally, to show the effectiveness of the method, we classify aCM semistable
bundles in rank $3$. 

\begin{prop}
Let $E$ be an aCM semistable normalized rank-$3$ bundle over $V_5$. Then $E$ is isomorphic to
one of the following
\begin{enumerate}[i)]
\item The dual universal $Q^*$ or the twisted universal $Q(-1)$;
\item $E^{(0,3)} \simeq \OO^{\oplus 3}$;
\item $E^{(1,2)} \simeq \OO \oplus E_L$;
\item $E^{(2,1)} \simeq \OO \oplus E_S$;
\item $E^{(3,0)} \simeq \Sym^2 U(-1)$ or a deformation of $\Sym^2 U(-1)$.
The deformation space of $\Sym^2 U(-1)$ has dimension $10$.
\end{enumerate}
\end{prop}

\begin{proof}
The first case is $c_1(E) = -1$. Then by
        $$0 \rt U^a \rt (Q^*)^b \rt E \rt 0$$
we get $3 = \rk(E) = 3b - 2a$, and $-1 = c_1(E) = -b+a$. Then $b=1$ and
$a=0$, that is $E = Q^*$. 
Clearly if $c_1(E)=-2$, then $c_1(E^*(-1)) = -1$, therefore $E =
Q(-1)$.

Now look at the case $c_1(E)=0$. Then we have
        $$0 \rt U^a \rt (Q^*)^b \oplus \OO^c \rt E \rt 0$$

Here we have $a=b$ and the rank $3$ says we have only the $4$ cases
        $$(a,c) \in \{(0,3),(1,2),(2,1),(3,0) \}$$

The first is the obvious $E \simeq \OO^3$. Cases $(1,2)$ and $(2,1)$ are described by
        $$\xymatrix@-2.5ex{
                                & {\OO} \ar[d] \ar@{=}[r]               & {\OO} \ar[d]          &
                        & & {\OO} \ar[d] \ar@{=}[r]             & {\OO} \ar[d]\\                
        U \ar@{=}[d] \ar[r]     & Q^* \oplus \OO^2 \ar[r] \ar[d]        & E^{(1,2)} \ar[d]        &
                        & U^2 \ar@{=}[d] \ar[r] & (Q^*)^2 \oplus \OO \ar[r] \ar[d]      & E^{(2,1)} \ar[d]        \\
        U \ar[r]                & Q^* \oplus \OO \ar[r]         & E_L                           &
                        & U^2 \ar[r]            & (Q^*)^2  \ar[r]                       & E_S
        }$$

So $E^{(1,2)} \simeq E_L \oplus \OO$ and $E^{(2,1)} \simeq E_S \oplus \OO$ since
$E_L$ and $E_S$ have no $\HH^1$. The last possibility is
        \begin{equation} \label{last}
        0 \rt U^3 \overset{\varphi}{\rt} (Q^*)^3 \rt E \rt 0
        \end{equation}

This corresponds to $\Sym ^2 U$. First normalize to $\Sym^2U(1)$, $c_1(\Sym^2 U(1))=0$. Then recall from Proposition
(\ref{UorO}) that $\Sym^2U$ is indeed aCM. Then using the dual of sequences (\ref{QQ}) and (\ref{universal}),
plus the symmetric square of (\ref{universal}), one computes 
        $$0 \rt Y_4 \rt Y_4 \oplus Y_2 \rt \HH^1(U \ts \Sym^2 U(1)) \rt 0$$

All the maps here are invariant so $\HH^1(U \ts \Sym^2 U(1)) \simeq
Y_2$. Furthermore, notice that 
        $$\HH^1(Q(-1) \ts \Sym^2 U(1)) \simeq  \HH^2(U(-2) \ts \Sym^2 U(2))
        \simeq \HH^1(U \ts \Sym^2 U(1))^* \simeq Y_2$$

Clearly $\Sym^2 U$ is stable. Now if we put the invariant map $\varphi _0
\in Y_2 \ts Y_2 \ts Y_2$ we get the resolution
        $$0 \rt U^3 \xrightarrow{\varphi_0} (Q^*)^3 \rt \Sym^2 U(1) \rt 0$$

Hence we generically get a stable aCM bundle as image of $\varphi$
with the same invariants as $\Sym ^2(U)(1)$. One may also notice that the copresentation of $E(-1)$ is
        $$0 \rt E(-1) \rt U^6 \rt (Q^*)^3 \rt 0$$

In the case of $\Sym^2 U$ the above arrow is an invariant element in
$Y_2 \ts Y_2 \ts (Y_4 \oplus \kk)$. By taking non invariant maps we choose a deformation of $\Sym^2U(-1)$.
The last statement follows, since tensorizing (\ref{last}) by $\Sym^2 U(-1)$
and taking cohomology yields $\hh^1(\End (\Sym^2 U)) =10$.
\end{proof}

\section*{Aknowledgements}
The author would like to express his gratitude to Enrique Arrondo, Laura Costa and Beatriz
Gra{\~{n}}a for many useful suggestions, and Giorgio Ottaviani for his
constant support. Also he would like to thank the Departments of Universitat de
Barcelona and Universidad Complutense de Madrid for the friendly hospitality. 

The author was partially supported by Italian Ministry funds
and the EAGER contract HPRN-CT-2000-00099.

\def\cprime{$'$} \def\cprime{$'$} \def\cprime{$'$} \def\cprime{$'$}
  \def\cprime{$'$}
\providecommand{\bysame}{\leavevmode\hbox to3em{\hrulefill}\thinspace}
\providecommand{\MR}{\relax\ifhmode\unskip\space\fi MR }
\providecommand{\MRhref}[2]{%
  \href{http://www.ams.org/mathscinet-getitem?mr=#1}{#2}
}
\providecommand{\href}[2]{#2}

\end{document}